\documentclass[12pt]{article}
\usepackage{latexsym,amsfonts,amssymb}
\usepackage[dvips]{color}
\usepackage{colordvi,multicol}
\usepackage{mathrsfs}    
\usepackage{amssymb}     
\usepackage{dutchcal}

\newtheorem{thm}{Theorem}[section]

\newtheorem{lem}[thm]{Lemma}
\newtheorem{pro}[thm]{Proposition}

\newtheorem{rem}[thm]{Remark}
\newtheorem{exa}[thm]{Example}

\begin{document}
\title{\bf  A lemma on a finite union-closed family of  finite sets and its applications}
\author{Ze-Chun Hu$^{1}$,  Yi-Ding Shi$^1$ and Qian-Qian Zhou$^{2,}$\footnote{Corresponding author: qianqzhou@yeah.net}\ \\ \\
 {\small $^1$ College of Mathematics, Sichuan  University,
 Chengdu 610065, China}\\
   {\small $^2$College of Science, Tianjin University of Technology, Tianjin 300384, China }}

\date{}
 \maketitle

\noindent{\bf Abstract}\quad
Suppose that $\mathscr{F}$ is a finite union-closed family of sets with $\cup_{A\in \mathscr{F}}A=\{1,2,\ldots,m\}$ and $m\geq 2$. Fix $i\in \{1,2,\ldots,m\}$ and denote
$\mathscr{G}:=\{A\backslash \{i\}: A\in \mathscr{F}\}$. For $j\in \{1,2,\ldots,m\}\backslash\{i\}$, let $\mathscr{G}_j:=\{A\in\mathscr{G}: j\in A\}$ and $\mathscr{F}_j:=\{A\in\mathscr{F}: j\in A\}$. In this note, we will prove a lemma which says that if $\frac{|\mathscr{G}_j|}{|\mathscr{G}|}\geq c\,(c\in (0,1])$, then
$\frac{|\mathscr{F}_j|}{|\mathscr{F}|}\geq \frac{1}{1+2(1-c)/c}$. Several applications of this lemma will be given.

\noindent{\bf Key words} The union-closed sets conjecture, Frankl's conjecture, Nagel's conjecture

\noindent{\bf Mathematics Subject Classification (2010)} 03E05, 05A05

\section{Introduction and main result}

In combinatorics, a very famous open problem is {\it the union-closed sets conjecture} (also called {\it Frankl's conjecture}, cf. \cite{Ri85, St86}), which says that  for any finite union-closed family of finite sets, other than the family consisting only of the empty set, there exists an element that belongs to at least half of the sets in the family.

A family $\mathscr{F}$ of sets is {\it union-closed} if, for any $A,B\in \mathscr{F}$, it holds that $A\cup  B\in \mathscr{F}$. For simplicity, denote $m=|\cup_{A\in \mathscr{F}}A|$ and $n=|\mathscr{F}|$. Hereafter, for any set $A$,  $|A|$  denotes the cardinal number of $A$.

Up to now, Frankl's conjecture is still open, and we only know it is true in some special cases. If a union-closed family $\mathscr{F}$ contains a set with one  or two elements,  Frankl's conjecture holds for $\mathscr{F}$ (\cite{SR89}). This result was extended by Poonen (\cite{Po92}). In addition,  the author in \cite{Po92}  proved that Frankl's conjecture holds if $m\leq 7$ or $n\leq 28$, and proved an equivalent lattice formulation of Frankl's conjecture. Bo\v{s}njak and Markovi\'{c} (\cite{BM08}) proved that Frankl's conjecture holds if $m\leq 11$.    Vu\v{c}kovi\'{c} and Zivkovi\'{c} (\cite{VZ17}) provided  a computer-assisted proof that Frankl's conjecture is true if  $m\leq 12$, which, together with Faro's result (\cite{Lo94-b}) (see also Roberts and Simpson \cite{RS10}), implies that Frankl's conjecture holds if  $n\leq 50$.  For further progress on Frankl's conjecture, we refer to  \cite{BS15}, \cite{EIL22}, \cite{Gi22}, \cite{HL20},  \cite{KPT23}, \cite{Ka17}, \cite{Li23},  \cite{Na23}, \cite{St21}, and  references therein.

Let $M_m=\{1,2,\ldots,m\}$ and $\mathscr{F}\subset 2^{M_m}=\{A: A\subset M_m\}$ with $\cup_{A\in\mathscr{F}}A=M_m$.  Suppose that $\mathscr{F}$ is union-closed. An obvious method to attack Frankl's conjecture is to use induction on the number $m$.
But it does not work, because for any $i\in \{1,2,\ldots,m\}$, the map: $A\in\mathscr{F} \mapsto A\backslash\{i\}$ is not injective in general.

Fix $i\in \{1,2,\ldots,m\}$. Define
\begin{eqnarray*}
\mathscr{G}:=\{A\backslash\{i\}: A\in \mathscr{F}\}.
\end{eqnarray*}
Suppose that $m\geq 2$ and $j\in \{1,2,\ldots,m\}\backslash\{i\}$. Denote
\begin{eqnarray*}
\mathscr{G}_j:=\{A\in \mathscr{G}: j\in A\},\ \mathscr{G}_{/ j}:=\{A\in \mathscr{G}: j\notin A\}.
\end{eqnarray*}
Similarly define $\mathscr{F}_j$ and $\mathscr{F}_{/ j}$.

The main {\it motivation} of this note is to consider the following question:

\quad\quad {\it What can we say about the ratio $\frac{|\mathscr{F}_j|}{|\mathscr{F}|}$ based on the information $\frac{|\mathscr{G}_j|}{|\mathscr{G}|}$?}

The main result is as follows:

\begin{lem}\label{lem}
If $\frac{|\mathscr{G}_j|}{|\mathscr{G}|}\geq c$ for some constant $c\in (0,1]$, then
$\frac{|\mathscr{F}_j|}{|\mathscr{F}|}\geq \frac{1}{1+2(1-c)/c}$.
\end{lem}

The basic idea for the proof of the above lemma is the equality (\ref{2.1}) below.
The inequality $\frac{|\mathscr{F}_j|}{|\mathscr{F}|}\geq \frac{1}{1+2(1-c)/c}$ in Lemma \ref{lem} is sharp, in the sense that if $\frac{|\mathscr{G}_j|}{|\mathscr{G}|}=c$ for some constant $c\in (0,1]$, then
$\frac{|\mathscr{F}_j|}{|\mathscr{F}|}=\frac{1}{1+2(1-c)/c}$ in some case. See the following example and Remark \ref{rem-2.2} below for the illustration.

\begin{exa}
Set $\mathscr{F}=\{\emptyset, \{1\}, \{3\}, \{1,2\},\{1,3\}, \{1,2,3\}\}.$
Take $i=1,j=2$. Then
$$
\mathscr{G}=\{\emptyset, \{2\}, \{3\}, \{2,3\}\},\ \mathscr{G}_2=\{\{2\}, \{2,3\}\},\ \mbox{and}\ \mathscr{F}_2=\{\{1,2\},\{1,2,3\}\}.
$$
Thus $c=\frac{|\mathscr{G}_2|}{|\mathscr{G}|}=\frac{1}{2}$ and $\frac{|\mathscr{F}_2|}{|\mathscr{F}|}=\frac{1}{3}=\frac{1}{1+2(1-c)/c}$.
\end{exa}

The rest of this note is organized as follows. The proof of Lemma \ref{lem} is provided  in Section 2 and some applications are presented in Section 3.

\section{Proof of Lemma \ref{lem}}\setcounter{equation}{0}

We will use the following elementary inequality. Its proof is obvious.

\begin{lem}\label{lem-2}
Suppose that $a,b,c,d$ are four positive numbers satisfying
$$
\frac{b}{a}\geq k,\ \mbox{and}\ \frac{d}{c}\geq k,
$$
for some positive constant $k$. Then $\frac{b+d}{a+c}\geq k$.
\end{lem}

\noindent {\bf Proof of Lemma \ref{lem}.} Define
\begin{eqnarray*}
&&x=|\{A\in \mathscr{F},i\notin A,j\in A\}\cap \{B\backslash\{i\}: B\in\mathscr{F}, i,j\in B\}|,\\
&&y=|\{A\in \mathscr{F},i\notin A,j\notin A\}\cap \{B\backslash\{i\}: B\in\mathscr{F}, i\in B,j\notin B\}|.
\end{eqnarray*}
Then $0\leq x\leq |\mathscr{G}_j|, 0\leq y\leq |\mathscr{G}_{/ j}|$, and
\begin{eqnarray}\label{2.1}
\frac{|\mathscr{F}_j|}{|\mathscr{F}|}=
\frac{|\mathscr{G}_j|+x}{|\mathscr{G}_j|+|\mathscr{G}_{/ j}|+x+y}.
\end{eqnarray}

By the assumption that $\frac{|\mathscr{G}_j|}{|\mathscr{G}|}\geq c$, i.e.,
$$
\frac{|\mathscr{G}_j|}{|\mathscr{G}_j|+|\mathscr{G}_{/  j}|}\geq c,
$$
we get that
\begin{eqnarray*}
|\mathscr{G}_{/ j}|\leq \frac{1-c}{c}|\mathscr{G}_j|,
\end{eqnarray*}
which together with the fact that $0\leq y\leq |\mathscr{G}_{/ j}|$ implies
\begin{eqnarray}\label{2.2}
\frac{|\mathscr{G}_j|}{|\mathscr{G}_j|+|\mathscr{G}_{/ j}|+y}\geq \frac{1}{1+2(1-c)/c}.
\end{eqnarray}
Without loss of generality, we can assume that $x>0$, then $\frac{x}{x}=1\geq \frac{1}{1+2(1-c)/c}$. Hence by (\ref{2.1}), (\ref{2.2}) and Lemma \ref{lem-2}, we obtain
$$\frac{|\mathscr{F}_j|}{|\mathscr{F}|}\geq \frac{1}{1+2(1-c)/c}.$$
The proof is complete.\hfill\fbox

\begin{rem}\label{rem-2.2}
By the above proof, we know that if $\frac{|\mathscr{G}_j|}{|\mathscr{G}|}= c\, (c\in (0,1])$,  $x=0$ and $y=|\mathscr{G}_{/ j}|$, then $\frac{|\mathscr{F}_j|}{|\mathscr{F}|}= \frac{1}{1+2(1-c)/c}$.

\end{rem}

\section{Applications}\setcounter{equation}{0}

In this section, we will give several applications of Lemma \ref{lem}.

\subsection{The equivalence of Frankl's conjecture and Nagel's conjecture}

By permuting the elements of the ground set $\{1,2,\ldots,m\}$, we can assume that
\begin{eqnarray}\label{3.1}
|\{F\in \mathscr{F}: 1\in F\}|\geq|\{F\in \mathscr{F}: 2\in F\}|\geq\cdots\geq |\{F\in \mathscr{F}: m\in F\}|.
\end{eqnarray}
Then  Frankl's conjecture says that
$$
|\{F\in \mathscr{F}: 1\in F\}|\geq \frac{1}{2}|\mathscr{F}|.
$$

Nagel \cite{Na23} introduced the following conjecture (we call it {\it Nagel's conjecture}):
for any $k\in \{1,2,\ldots,m\}$, it holds that
\begin{eqnarray}\label{3.2}
|\{F\in \mathscr{F}: k\in F\}|\geq \frac{1}{2^{k-1}+1}|\mathscr{F}|.
\end{eqnarray}

Obviously,  Nagel's conjecture implies Frankl's conjecture. Das and Wu \cite{DS24} proved that Nagel's conjecture holds for $k\geq 3$ and for $k=2$ under additional conditions.

By Lemma \ref{lem}, we have the following result.

\begin{pro}  Frankl's conjecture is equivalent to Nagel's conjecture.
\end{pro}
{\bf Proof.} It is sufficient to demonstrate  the necessity. Suppose that Frankl's conjecture is true. Then, for the $\mathscr{F}$ above, we have
$$
|\{F\in \mathscr{F}: 1\in F\}|\geq \frac{1}{2}|\mathscr{F}|.
$$
Without loss of generality, we assume that $m\geq 2$. Define
$$
\mathscr{G}:=\{A\backslash \{1\}: A\in \mathscr{F}\}.
$$
 Subsequently, $\mathscr{G}$ is a union-closed family of sets satisfying $\cup_{A\in \mathscr{G}}A=\{2,3,\ldots,m\}$. By the assumption that Frankl's conjecture is true, we know that there exists $i\in \{2,3,\ldots,m\}$ such that
$$
|\{A\in \mathscr{G}: i\in A\}|\geq \frac{1}{2}|\mathscr{G}|.
$$
Then by Lemma \ref{lem}, we get that
$$
|\{A\in \mathscr{F}: i\in A\}|\geq \frac{1}{1+\frac{2(1-1/2)}{1/2}}|\mathscr{F}|=\frac{1}{3}|\mathscr{F}|
=\frac{1}{2^{2-1}+1}|\mathscr{F}|,
$$
which together with (\ref{3.1}) implies that
\begin{eqnarray*}
|\{F\in \mathscr{F}: 2\in F\}|\geq \frac{1}{2^{2-1}+1}|\mathscr{F}|,
\end{eqnarray*}
i.e., (\ref{3.2}) holds for $k=2$.

Notice that
$$
\frac{1}{1+\frac{2(1-\frac{1}{2^{k-1}+1})}{\frac{1}{2^{k-1}+1}}}=\frac{1}{2^k+1}
=\frac{1}{2^{(k+1)-1}+1}.
$$
Then, using  the reduction method, we can obtain  that for any $k=3,\ldots,m$, (\ref{3.2}) holds. Thus,  Nagel's conjecture is valid. The proof is complete. \hfill\fbox

\subsection{A complement to  \cite[Lemma 2.4]{Na23}}

First, we recall \cite[Lemma 2.4]{Na23} as follows:

\begin{lem}\label{lem-3.3}
 (Nagel)  For any $A\in \mathscr{F}$ with $|A|\geq 1$, and any $x\in A$, it holds that
$$
|\{F\in\mathscr{F}: x\in F\}|\geq \frac{1}{2^{|A|-1}+1}|\mathscr{F}|.
$$
\end{lem}

If $|A|=1$, then $\frac{1}{2^{|A|-1}+1}=\frac{1}{2}$ and thus the inequality in Lemma \ref{lem-3.3} is the best in this case. If $|A|\geq 2$, we obtain  the following result which can be regarded as a complement to Lemma \ref{lem-3.3}.

\begin{pro}\label{pro-3.4}
 For any $A\in \mathscr{F}$ with $|A|\geq 2$, there exists $y\in A$ such that
\begin{eqnarray}\label{pro-3.4-a}
|\{F\in\mathscr{F}: y\in F\}|\geq \frac{1}{2^{|A|-2}+1}|\mathscr{F}|.
\end{eqnarray}
\end{pro}
{\bf Proof.} If $|A|=2$, it is well known that one element $y\in A$ exists   such that $$
|\{F\in \mathscr{F}: y\in F\}|\geq \frac{1}{2}|\mathscr{F}|,
$$
i.e., (\ref{pro-3.4-a}) holds in this case.

If $|A|=3$, we express $A$ as $\{x_1,x_2,x_3\}$. Define
$$
\mathscr{G}:=\{B\backslash \{x_1\}: B\in \mathscr{F}\}.
$$
Then $\mathscr{G}$ is a finite union-closed family with $\cup_{B\in \mathscr{G}}B=\{1,2,\cdots,m\}\backslash\{x_1\}$ and $\{x_2,x_3\}\in \mathscr{G}$. Thus,  $k\in \{2,3\}$  exists  such that
$$
|\{B\in \mathscr{G}: x_k\in B\}|\geq \frac{1}{2}|\mathscr{G}|.
$$
Then by Lemma \ref{lem}, we get that
$$
|\{B\in \mathscr{F}: x_k\in B\}|\geq \frac{1}{1+2(1-\frac{1}{2})/\frac{1}{2}}|\mathscr{F}|=\frac{1}{3}|\mathscr{F}|
=\frac{1}{2^{3-2}+1}|\mathscr{F}|,
$$
i.e., (\ref{pro-3.4-a}) holds in this case.

For $|A|\geq 4$, by the reduction method, we can easily obtain the result.\hfill\fbox

\subsection{About one stronger version of Frankl's conjecture}

For any $k=1,2,\ldots,m$, denote $\mathscr{M}_k=\{A\in 2^{M_m}: |A|=k\}.$ Define
$$
T(\mathscr{F})=\inf\{1\leq k\leq m: \mathscr{F}\cap \mathscr{M}_k\neq \emptyset\}.
$$
Then $1\leq T(\mathscr{F})\leq m$.  By virtue of $T(\mathscr{F})$, Cui and Hu \cite{CH21} introduced two  stronger versions of Frankl's conjecture, one of which is as follows:

{\bf $S_2$-version:} if  $m\geq 2$ and $T(\mathscr{F})\geq 2$, then there exist at least two elements in $M_m$ that belong to at least half of the sets in $\mathscr{F}$.

  It is easy to know that $S_2$-version implies  Frankl's conjecture. Similar to the result in  Gilmer \cite{Gi22},  we can ask the following question:

{\it Question 1.} If $m\geq 2$, are there two positive constants $c_1,c_2$ with $c_1\geq c_2$ such that there exist two different elements $i,j\in \{1,2,\cdots,m\}$ such that
$$
\frac{|\mathscr{F}_i|}{|\mathscr{F}|}\geq c_1,\ \mbox{and}\ \frac{|\mathscr{F}_j|}{|\mathscr{F}|}\geq c_2?
$$

\begin{rem} As to the constants $c_1$ and $c_2$ in Question 1, we have the following conjecture and results.

(i) By virtue of Subsection 3.1, we conjecture  that the best possible general result is that $c_1=\frac{1}{2}, c_2=\frac{1}{3}$.

(ii) If $T(\mathscr{F})=2$, then by Lemma \ref{lem} or Proposition \ref{pro-3.4}, we can take  $c_1=\frac{1}{2}, c_2=\frac{1}{3}$.

(iii) S-Frankl's conjecture is equivalent to that if $T(\mathscr{F})\geq 2$, then $c_1=c_2=\frac{1}{2}$.

(iv) From  the  results in Liu \cite{Li23}, we can take $c_1=0.38234$. Then by Lemma \ref{lem}, we can take $c_2=\frac{1}{1+2(1-0.38234)/0.38234}\approx 0.23635$.

\end{rem}

\vskip 0.5cm
{ \noindent {\bf\large Acknowledgments}

 This study  was supported by the National Natural Science Foundation of China (Grant Nos. 12171335  and, 12301603).

\bigskip

\noindent {\bf Author Contributions}: All authors contribute equally.

\noindent {\bf Conflict of Interest}: The authors declare no competing interests.

\end{document}